\documentclass{article}
\usepackage{graphicx} 
\usepackage[T1]{fontenc}
\usepackage{amsfonts}
\usepackage{csquotes}
\usepackage[english]{babel}
\usepackage[a4paper,top=2cm,bottom=2cm,left=2cm,right=2cm,marginparwidth=1.75cm]{geometry}
\usepackage{amsmath}
\usepackage{tikz}
\usetikzlibrary{positioning}
\usepackage{mathtools}
\usepackage[font=small,skip=7pt]{caption}
\usepackage[aboveskip=-2pt]{subcaption}
\usepackage{float}
\usepackage[colorlinks=true, allcolors=blue]{hyperref}
\usepackage{comment}
\usepackage{natbib}

\DeclareMathOperator{\Tr}{Tr}

\renewcommand{\i}{\mathrm{i}}

\title{Hopf bifurcation and periodic solutions in a sustainable supply chain defined by a planar system of ordinary differential equations}
\author{Alessia Andò and Dimitri Breda}

\begin{document}
\maketitle

\begin{abstract}
We analyze the stability of the equilibria and bifurcations of a planar system of ordinary differential equations describing the product-resource interaction in a sustainable supply chain. While periodic behavior in supply chain models has often been documented either in systems of higher dimensions or in delayed systems, here a Hopf bifurcation arises in a planar and delay-free system. We show that the interior equilibrium loses stability through a supercritical Hopf bifurcation as the environmental capacity exceeds a critical threshold that depends on the maximum production rate, the resource level at which production reaches half its maximum rate, and the demand and remanufacturing rates. We derive this threshold explicitly and provide estimates for the amplitude and period of the periodic solution close to the Hopf bifurcation, by means of the resulting Hopf normal form. We then further substantiate our analytical results through numerical simulations.
\end{abstract}

\noindent
\textbf{Keywords:} supply chain modeling, sustainable supply chains, periodic solution, Hopf bifurcations, normal form theory, environmental capacity

\smallskip
\noindent{\bf{2020 Mathematics Subject Classification:}} 34C20, 34C23, 34C25, 34C60, 90B30, 91B76

\section{Introduction}
Supply chains are among the most consequential engineered systems of the contemporary economy, in that they coordinate the flow of materials and products across companies and borders. Their performance affects cost, ability to meet the customers' demand and, more importantly, sustainability. Indeed, amid the increasing awareness of resource efficiency, waste reduction and environmental cost of production \citep{ONU}, supply chains have become a primary factor in whether sustainability targets are met in practice.

Being feedback-driven systems, supply chains have long been studied using tools from nonlinear dynamics. Ordering policies, production rates and resource exploitation can be modeled as possible components of a dynamical system whose equilibria, stability, and bifurcations can be analyzed with standard methods for differential equations  \citep{Forrester1961,Sterman1989}. Determining when such systems settle into a stable equilibrium or limit cycle is not only relevant to economy, but also to the efficient use of natural resources.

Understanding and controlling oscillations in these systems is a persistent challenge. The so-called \emph{bullwhip effect} \citep{Lee1997} explains how information distortion at different levels can be the cause of oscillations, and has a sustainability dimension: overproduction leads to increased resource consumption and carbon emission. Strongly related is the concept of \emph{green bullwhip effect} \citep{Lee2014} through which environmental requirements upstream in the same way as demand variability does. This motivates the interest in understanding the birth of oscillations in supply chains, not only due to the willingness to protect profitability, but above all in order to manage waste and emissions associated with erratic production and transport.

An equilibrium that loses stability as a parameter crosses a critical value can give rise, through a Hopf bifurcation, to a self-sustained periodic orbit \citep{Guckenheimer2013}. A substantial body of recent work has explored this mechanism in supply chain models, but most of it uses deliberately systems of dimension 3 or more, in order to exhibit chaos in addition to periodicity \citep{Anne2009,Mondal2019,Zheng2022,Johansyah2024a,Johansyah2024b,Johansyah2025,Shi2026}. Another strand of the literature considers the incorporation of delays, which have different interpretations depending on the context. For instance, in supply chain financing models, the delay represents the time required for information transmission and strategic decisions between firms and financial institutions. These delayed feedback mechanisms have been shown to destabilize equilibria and induce Hopf bifurcations \citep{Chen2020}. Oscillations in supply chain models may also stem from delays in manufacturers' and retailers' quality and service decisions, where control actions are based on historical market information, such as past prices and consumer demand \citep{Han2021}. Sufficiently large information-processing and verification lags between banks and B2B platforms
may likewise trigger Hopf bifurcations and oscillatory behavior in the strategic dynamics \citep{Yan2026}.

What is comparatively scarce in this literature, to the best of our knowledge, is the study of delay-free two-dimensional ordinary differential equations (ODEs) of supply chain models in which a Hopf bifurcation arises from the nonlinear interaction between the state variables. The present paper aims to address this gap.

The model we analyze is cast in terms of a product variable and a resource variable, and was originally proposed by Wang and Gunasekaran, who introduced a nonlinear dynamical system to describe the interaction between supply chains and their environmental impact, and used it to study how design production capacity, environmental cost, and demand rate affect performance \citep{Wang2015}. As originally formulated, the model has three state variables, apparently in line with the three-dimensional constructions discussed above; however, one of these variables is uncoupled from the other two, so that the dynamics relevant to the interaction between products and resources reduce to a two-dimensional autonomous ODE. That paper established existence and uniqueness of solutions and carried out an equilibrium and local stability analysis; however, the possibility of sustained periodic solutions was not identified, so the oscillatory regime documented in the present paper had, to our knowledge, gone undetected. Our aim here is to revisit this stability analysis: after identifying the equilibria, we derive the conditions on the model parameters under which a Hopf bifurcation occurs, and establish the direction and stability of the resulting periodic orbit via the first Lyapunov coefficient. Then, by means of the Hopf normal form, we obtain an estimate of the amplitude and period close to the Hopf bifurcation, and substantiate our findings through numerical simulations. In particular, we fiund that increasing the environmental capacity or the maximum production rate destabilizes the equilibrium and pushes the system toward oscillation, while increasing the resource level at which production reaches half its maximum rate, or increasing the remanufacturing rate, restores stability.

The remainder of the paper is organized as follows.  In Section \ref{sec:model} we introduce the supply chain model that we analyze in the rest of the paper. Then, we find its equilibria and study their stability in Section \ref{sec:equilibria}. In Section \ref{sec:hopf} we focus on the nontrivial equilibrium and conclude that its stability change corresponds to a supercritical Hopf bifurcation. We then obtain an estimate for the amplitude and period of the periodic solution slightly past the Hopf bifurcation; the most technical computations are collected in Appendix \ref{sec:c2}. In Section \ref{sec:simulations} we provide numerical evidence of all the results obtained in Section \ref{sec:hopf}. We conclude in Section \ref{sec:conclusion} with some final remarks about future extension of the present work.
\section{Model formulation}\label{sec:model}
The model proposed by \citet{Wang2015}, and then extended with a stochastic component by \citet{DAmbrosio2026}, adapts concepts from microbial biology to track the interaction between supply chain operations and the natural environment. It involves three (operational, ecological, and financial state) variables, by means of the ODE
\begin{equation*}
\left\{
\begin{aligned}
x_1' &= \left(\bar{p}\,\frac{x_2}{x_2+K_P}-(d-m)\right)x_1, \\[6pt]
x_2' &= bx_2\!\left(1-\frac{x_2}{K_E}\right)-\frac{ax_1 x_2}{x_2+K_P}\\[6pt]
x_3'&=\frac{(r-c)x_1 x_2}{x_2+K_P}-fx_1-\mu x_2,
\end{aligned}
\right.
\end{equation*}
where $x_1$ represents the total amount of products inside the supply chain network, $x_2$ the total units of available natural resources provided by the environment, and $x_3$ the accumulated financial profit of the supply chain. Concerning the production parameters, appearing in the first equation,
\begin{itemize}
    \item $\bar{p}$ is the maximum production rate that the supply chain is theoretically able to reach, provided abundant resources, and is an upper bound to the effective production rate $\bar{p}\,\frac{x_2}{x_2+K_P}$. The latter formula is an adaptation of Monod kinetics, now widely used to model biomass growth in bioreactors \citep{Liu2017,Monod1949};
    \item $K_P$ represents the resource level at which the factories run at exactly half of their maximum design capability ($\bar{p}/2$);
    \item $d$ is the demand rate, at which customers buy or pull finished products out of the supply chain system;
    \item $m$ is the remanufacturing rate, at which damaged or obsolete products are collected and then recycled.
\end{itemize}
As for the resource parameters in the second equation,
\begin{itemize}
    \item $b$ represents the natural regeneration/birth rate of the raw environmental resources;
    \item $K_E$ represents the environmental carrying capacity, beyond which a biological population or resource cannot naturally regenerate due to structural constraints, and defines a logistic growth mechanism;
    \item $a$ is the maximum speed or intensity at which the active manufacturing infrastructure consumes raw materials to sustain production.
\end{itemize}
Finally, the parameters involved in the financial profit (third equation) have the following meaning:
\begin{itemize}
    \item $r$ is the revenue generated by the supply chain for every single unit of product sold;
    \item $c$ is the basic operational cost required to manufacture, process and distribute one unit of product;
    \item $f$ is the financial penalty incurred per unit of resource usage. This can be imagined as the result of carbon taxes, ecological rehabilitation expenses and possible regulatory fines;
    \item $\mu$ is the cost to collect, disassemble, process and recycle each end-of-life unit of product.
\end{itemize}
Due to the meaning of the parameters, it is assumed that all of them are positive, and that $d>m$. Observe that the financial profit equation serves crucial conceptual purposes in \citep{Wang2015}. Indeed, the authors mention that one of the contributions is to facilitate the optimal control in supply chain operations, where $x_3$ would be the variable to optimize. However, since from a structural point of view the equation is uncoupled from the other two, it is sufficient for us to carry out our equilibria and stability analysis of the reduced system
\begin{equation}\label{eq:Wang}
\left\{
\begin{aligned}
x_1' &= \left(\bar{p}\,\frac{x_2}{x_2+K_P}-(d-m)\right)x_1, \\[6pt]
x_2' &= bx_2\!\left(1-\frac{x_2}{K_E}\right)-\frac{ax_1 x_2}{x_2+K_P}.
\end{aligned}
\right.
\end{equation}
We define \[
  p^* := \frac{\bar{p}}{d-m}.
\]
\section{Analysis of equilibria}\label{sec:equilibria}
In this section, we compute the coordinates of the equilibria of \eqref{eq:Wang} and then study their stability by means of the principle of linearzied stability. In order to obtain the coordinates of the equilibria, we compute the nullclines.
\paragraph{$x_1$-nullclines.} Either $x_1=0$, or (dividing by $x_1$)
\begin{equation*}
  \bar{p}\,\frac{x_2}{x_2+K_P} = d-m
  \iff
  \frac{x_2}{x_2+K_P} = \frac{1}{p^*}
  \iff
  x_2(p^*-1) = K_P.
\end{equation*}
\paragraph{$x_2$-nullclines.} Either $x_2=0$, or (dividing by $x_2$)
\begin{equation*}
  x_1 = \frac{b(x_2+K_P)}{a}\left(1-\frac{x_2}{K_E}\right).
\end{equation*}
Thus, we have the trivial equilibrium $E_0:=(0,0)$ and the equilibrium in the absence of production $E_1:=(0,K_E)$ which always exist, as well as $E_2:=(\overline{x}_1,\overline{x}_2)$ for
\begin{equation}\label{eq:E2}
  \overline{x}_1 := \frac{bK_P p^*}{a(p^*-1)}\!\left(1-\frac{K_P}{K_E(p^*-1)}\right) > 0, \qquad
  \overline{x}_2 := \frac{K_P}{p^*-1}>0.
\end{equation}
$E_2$ exists if and only if
\begin{equation*}
p^*>1
\end{equation*}
and 
\[K_E(p^*-1)>K_P.\]
\subsection{Stability of equilibria}
The Jacobian of~\eqref{eq:Wang} is
\begin{equation*}
J(x_1,x_2)=
\begin{pmatrix}
\bar{p}\,\dfrac{x_2}{x_2+K_P}-(d-m)
& \bar{p}\,\dfrac{K_P}{(x_2+K_P)^2}\,x_1 \\[10pt]
-\dfrac{a x_2}{x_2+K_P}
& b\left(1-\dfrac{2x_2}{K_E}\right)-\dfrac{a x_1 K_P}{(x_2+K_P)^2}
\end{pmatrix}.
\end{equation*}
As for the trivial equilibrium,
\[
J(E_0)=
\begin{pmatrix}
-(d-m) & 0 \\
0 & b
\end{pmatrix},
\]
with eigenvalues $\lambda_1=-(d-m)<0$ and $\lambda_2=b>0$. So, $E_0$ is always a saddle (unstable). As for the equilibrium in the absence of production,
\[
J(E_1)=
\begin{pmatrix}
\bar{p}\,\dfrac{K_E}{K_E+K_P}-(d-m) & 0 \\[10pt]
-\dfrac{aK_E}{K_E+K_P} & -b
\end{pmatrix}
\]
has eigenvalues
\[
\lambda_1=\bar{p}\,\frac{K_E}{K_E+K_P}-(d-m)=(d-m)\left(\frac{p^*K_E}{K_E+K_P}-1\right),\qquad\lambda_2 =-b<0.
\]
Thus, $E_1$ is (locally asymptotically) stable if and only if
\[\frac{p^*K_E}{K_E+K_P}<1\iff p^*K_E < K_E+K_P\iff K_E(p^*-1)<K_P.\]
Thus, if $p^*\leq 1$, then $E_1$ is stable for all values of $K_E,\,K_P$, while $E_2$ does not exist. If $p^*>1$, $E_1$ is stable $\iff E_2$ does not exist. This means that at $K_E(p^*-1)=K_P$ there is a transcritical bifurcation, where $E_2$ is born and $E_1$ loses stability. As for the equilibrium $E_2$, defined by \eqref{eq:E2},
\begin{align*}
J(E_2)=&
\begin{pmatrix}
0 & \bar{p}\,\dfrac{K_P}{(\overline{x}_2+K_P)^2}\,\overline{x}_1 \\[4pt]
-\dfrac{a\overline{x}_2}{\overline{x}_2+K_P} & b - \dfrac{2b\overline{x}_2}{K_E} - \dfrac{a\overline{x}_1 K_P}{(\overline{x}_2+K_P)^2}
\end{pmatrix}=\begin{pmatrix}
0 & \bar{p}\,\dfrac{bK_P}{a(\overline{x}_2+K_P)}\left(1-\dfrac{\overline x_2}{K_E}\right) \\[4pt]
-\dfrac{a}{p^*} &  b\left[1 - \dfrac{2\overline{x}_2}{K_E}
     - \dfrac{K_P\!}{\overline{x}_2+K_P}\left(1-\dfrac{\overline{x}_2}{K_E}\right)\right]
\end{pmatrix}\\[1em] \notag
=&\begin{pmatrix}
0 & \dfrac{(d-m)b(p^*-1)}{a}\left(1-\dfrac{K_P}{K_E(p^*-1)}\right) \\[4pt]
-\dfrac{a}{p^*} &  b\left[1 - \dfrac{2K_P}{K_E(p^*-1)}
     - \dfrac{p^*-1}{p^*}\left(1-\dfrac{K_P}{K_E(p^*-1)}\right)\right]
\end{pmatrix}&,
\end{align*}
where the equalities follow from \eqref{eq:E2}. Its determinant is 
\begin{equation}\label{eq:detJE2}
\Delta:=\det(J(E_2))=\dfrac{(d-m)b(p^*-1)}{p^*}\left(1-\dfrac{K_P}{K_E(p^*-1)}\right)>0
\end{equation}
for all possible parameter values.
The bottom right entry (and, thus, the trace) of $J(E_2)$ can be further simplified as
\begin{align*}
b\left[\frac{1}{p^*}
     + \frac{K_P}{K_E}\!\left(\frac{1}{p^*}-\frac{2}{p^*-1}\right)\right]
  = b\left[\frac{1}{p^*}
     - \frac{K_P(p^*+1)}{K_E p^*(p^*-1)}\right] 
  = \frac{b}{p^*K_E}\left(K_E - K_P\frac{p^*+1}{p^*-1}\right).\end{align*}
Thus, $E_2$ exists and is stable $\iff K_P\dfrac{1}{p^*-1}<K_E<K_P\dfrac{p^*+1}{p^*-1}$. The characteristic polynomial of $J(E_2)$ is
\begin{equation}\label{eq:charODE}
  \lambda^2 - \frac{b}{p^*K_E}\left(K_E - K_P\frac{p^*+1}{p^*-1}\right)\lambda + \dfrac{(d-m)b(p^*-1)}{p^*}\left[1-\dfrac{K_P}{K_E(p^*-1)}\right] = 0,
\end{equation}
The eigenvalues at $K_E=K_E^H:=K_P\dfrac{p^*+1}{p^*-1}$ are imaginary,  in particular
\begin{align*}\lambda_{1,2}^H=\pm\i\sqrt{\dfrac{(d-m)b(p^*-1)}{p^*}\left(1-\dfrac{K_P}{K_E(p^*-1)}\right)}=\pm\i\sqrt{\dfrac{(d-m)b(p^*-1)}{p^*}\left(1-\dfrac{1}{p^*+1}\right)}=\pm\i\sqrt{b(d-m)\dfrac{p^*-1}{p^\ast+1}}.
\end{align*}
This indicates the presence of a Hopf bifurcation at $K_E^H$. Indeed, by defining
\begin{equation}\label{eq:alpha}
\alpha=\alpha(K_E):=\frac{\Tr(J(E_2))}{2}=\frac{b}{2p^*K_E}\left(K_E - K_P\frac{p^*+1}{p^*-1}\right),
\end{equation}
we have $\alpha(K_E^H)=0$.
\section{Hopf bifurcation analysis}\label{sec:hopf}
In this section we analyze the Hopf bifurcation just found. In particular, we prove that it is supercritical, and then obtain an estimate for the amplitude and period of the periodic solution close to the Hopf bifurcation. 
The roots of \eqref{eq:charODE} are
\begin{align*}\lambda_{1,2}=&\alpha\pm\i\beta, \quad\beta:=\sqrt{-\alpha^2+\Delta},
\end{align*}
for $\Delta$ as in \eqref{eq:detJE2} and $\alpha$ as in \eqref{eq:alpha}. Slightly past the Hopf bifurcation, i.e., for $K_E>K_E^H$, we have $\alpha>0$. Thus, we can write 
\[
J(E_2)=
\begin{pmatrix}0 & \dfrac{(\alpha^2+\beta^2)p^*}{a}\\[8pt]
        -\dfrac{a}{p^*} & 2\alpha\end{pmatrix}.
\]
To obtain the Hopf normal form, we first translate $E_2$ to the origin via the change of coordinates $u = x_1-\overline{x}_1$, $v = x_2-\overline{x}_2$, obtaining
\[\left\{
\begin{aligned}
u' &= F(u,v):=\left(\bar{p}\,\frac{v+\overline x_2}{v+\overline x_2+K_P}-(d-m)\right)(u+\overline x_1), \\
v' &= G(u,v):=b(v+\overline x_2)\!\left(1-\frac{v+\overline  x_2}{K_E}\right)-\frac{a(u+\overline x_1)(v+\overline x_2)}{v+\overline x_2+K_P}.
\end{aligned}
\right.\]
Thus, the Taylor expansion up to order 3 reads
\begin{equation}\label{eq:uv3}
\left\{\begin{aligned}
u' &=\dfrac{(\alpha^2+\beta^2)p^*}{a}v+ F_{uv}uv + \dfrac{1}{2}F_{vv}v^2
  + \dfrac{1}{2}F_{uvv}uv^2 + \dfrac{1}{6}F_{vvv}v^3,\\[4pt]
  v'&=-\dfrac{a}{p^*}u+ 2\alpha v+ G_{uv}uv + \dfrac{1}{2}G_{vv}v^2
  + \dfrac{1}{2}G_{uvv}uv^2 + \dfrac{1}{6}G_{vvv}v^3,
  \end{aligned}\right.
  \end{equation}
  where, at $\alpha=0$ (i.e., for $K_E=K_E^H$), the coefficients are
  \[
  \begin{array}{rl}
F_{uv}=&\dfrac{\overline p}{K_P}\left(1-\dfrac{1}{p^*}\right)^2,\quad F_{vv}=-\dfrac{2\overline pb}{aK_P}\dfrac{(p^*-1)^2}{p^*(p^*+1)},\quad F_{uvv}=-\dfrac{2\overline p}{K_P^2}\left(1-\dfrac{1}{p^*}\right)^3,\quad F_{vvv}=\dfrac{6\overline pb}{aK_P^2}\dfrac{(p^*-1)^3}{p^{*2}(p^*+1)},\\
G_{uv}=&-\dfrac{a}{K_P}\left(1-\dfrac{1}{p^*}\right)^2,\quad G_{vv}=-\dfrac{2b(p^*-1)}{K_Pp^*(p^*+1)},\quad G_{uvv}=-\dfrac{2a}{K_P^2}\left(1-\dfrac{1}{p^*}\right)^3,\quad G_{vvv}=-\dfrac{6b(p^*-1)^3}{p^{*2}(p^*+1)K_P^2}.
  \end{array}
  \]
An eigenvector of $J(E_2)$ for $\lambda_1=\alpha+\i\beta$ is given by \[
v=\begin{pmatrix}\dfrac{(\alpha^2+\beta^2)p^*}{a}\\\lambda_1\end{pmatrix}=\begin{pmatrix}\dfrac{(\alpha^2+\beta^2)p^*}{a}\\\alpha\end{pmatrix}+\i\begin{pmatrix}0\\\beta\end{pmatrix}.\]
Thus, in order to translate coordinates to get the system into the canonical form for the Hopf theorem, we use the change of coordinates matrix
\[
T= \begin{pmatrix} \dfrac{(\alpha^2+\beta^2)p^*}{a} &0\\
\alpha &\beta\end{pmatrix}=: \begin{pmatrix} P &0\\
\alpha &\beta\end{pmatrix},\quad\text{with}\quad T^{-1}=\begin{pmatrix} P^{-1} &0\\
-\alpha P^{-1}\beta^{-1} & \beta^{-1}\end{pmatrix}\quad\text{and}\quad D=\begin{pmatrix} \alpha &\beta\\
-\beta &\alpha\end{pmatrix},\]
such that
\[
T^{-1}J(E_2)T= \begin{pmatrix} 0 &1\\
-\dfrac{a\beta^{-1}}{p} &\alpha\beta^{-1}\end{pmatrix}T=D.
\]
Now our new coordinates $(\xi,\eta)$ are defined by
\[\begin{pmatrix}u\\v\end{pmatrix}=T\begin{pmatrix}\xi\\\eta\end{pmatrix}=\begin{pmatrix}P\xi\\\alpha\xi+\beta\eta\end{pmatrix}.\]
Thus,
\begin{align*}
\begin{pmatrix}\xi'\\\eta'\end{pmatrix}=&T^{-1}\begin{pmatrix}u'\\v'\end{pmatrix}
=\begin{pmatrix}\alpha\xi+\beta\eta+f(\xi,\eta)\\
-\beta\xi+\alpha\eta +g(\xi,\eta)\end{pmatrix}
\end{align*}
where, at $\alpha=0$, we have
\[
\begin{aligned}
f(\xi,\eta)=&\beta F_{uv}\xi\eta+\dfrac{\beta^2}{2P}F_{vv}\eta^2+\dfrac{\beta^2}{2}F_{uvv}\xi\eta^2+\dfrac{\beta^3}{6P}F_{vvv}\eta^3+\cdots\\
g(\xi,\eta)=&PG_{uv}\xi\eta+\dfrac{\beta}{2}G_{vv}\eta^2+\dfrac{P\beta}{2}G_{uvv}\xi\eta^2+\dfrac{\beta^2}{6}G_{vvv}\eta^3+\cdots
\end{aligned}
\]
Thus, from \cite[Equation (3.4.11)]{Guckenheimer2013}, the first Lyapunov coefficient is
\begin{equation}\label{eq:lyap1}
\setlength\arraycolsep{0.1em}\begin{array}{rcl}
  l_1 &=& \cfrac{1}{16}\left((f_{\xi\eta\eta}+g_{\eta\eta\eta})
  +\cfrac{1}{\beta}(g_{\xi\eta}g_{\eta\eta}
  -f_{\eta\eta}(f_{\xi\eta}+g_{\eta\eta}))\right)\\[4mm]
  &=&\cfrac{1}{16}\left(\beta^2(F_{uvv}+G_{vvv})
  +\cfrac{1}{\beta}\left(P\beta G_{uv}G_{vv}
  -\cfrac{\beta^2}{P}F_{vv}(\beta F_{uv}+\beta G_{vv})\right)\right)\\[4mm]
  &=&\cfrac{1}{16}\Bigg(\beta^2\left(-\cfrac{2(p^*-1)^3}{K_P^2p^{*2}}
\left((d-m)+\cfrac{3b}{p^*+1}\right)\right)\\[4mm]
  &&+\dfrac{2Pab(p^*-1)^3}{K_P^2p^{*3}(p^*+1)}
  -\cfrac{2\beta^2b(d-m)(p^*-1)^3}{Pap^*(p^*+1)K_P^2}\left((d-m)(p^*-1)-\cfrac{2b}{p^*+1}\right)\Bigg)\\[4mm]
   &=&\cfrac{1}{16}\Bigg(\beta^2\left(-\cfrac{2(p^*-1)^3}{K_P^2p^{*2}}
\left((d-m)+\cfrac{3b}{p^*+1}\right)\right)\\[4mm]
  &&+\dfrac{2\beta^2b(p^*-1)^3}{K_P^2p^{*2}(p^*+1)}
  +\cfrac{2b(d-m)(p^*-1)^3}{p^{*2}(p^*+1)K_P^2}\left((d-m)(p^*-1)-\cfrac{2b}{p^*+1}\right)\Bigg)\\[4mm]
    &=&\cfrac{(p^*-1)^3}{8K_P^2p^{*2}}\left(\beta^2\left(-(d-m)-\cfrac{3b}{p^*+1}\right)
  +\dfrac{\beta^2b}{p^*+1}
  +\cfrac{b(d-m)}{(p^*+1)}\left((d-m)(p^*-1)-\cfrac{2b}{p^*+1}\right)\right)\\[4mm]
  &=&\cfrac{(p^*-1)^3}{8K_P^2p^{*2}}\left(-\beta^2\left(d-m+\cfrac{2b}{p^*+1}\right)
  +\cfrac{b(d-m)}{(p^*+1)}\left((d-m)(p^*-1)-\cfrac{2b}{p^*+1}\right)\right)\\[4mm]
    &=&\cfrac{(p^*-1)^3}{8K_P^2p^{*2}}\left(-\cfrac{b(d-m)(p^*-1)}{p^*+1}\left(d-m+\cfrac{2b}{p^*+1}\right)
  +\cfrac{b(d-m)}{(p^*+1)}\left((d-m)(p^*-1)-\cfrac{2b}{p^*+1}\right)\right)\\[4mm]
      &=&\cfrac{(p^*-1)^3}{8K_P^2p^{*2}}\left(-\cfrac{b(d-m)(p^*-1)}{p^*+1}\left(\cfrac{2b}{p^*+1}\right)
  +\cfrac{b(d-m)}{(p^*+1)}\left(-\cfrac{2b}{p^*+1}\right)\right)\\[4mm]
        &=&\cfrac{b^2(d-m)(p^*-1)^3}{4K_P^2p^{*2}(p^*+1)^2}\left(-(p^*-1)
  -1\right)\\[4mm]
  &=&-\cfrac{b^2(d-m)(p^*-1)^3}{4K_P^2p^{*}(p^*+1)^2}<0,
\end{array}
\end{equation}
  confirming that the bifurcation is supercritical. Note that, if $K_E=(1+h)K_E^H$ for some $h>0$ sufficiently small, then
  \[
  \alpha(K_E)=\frac{b}{2p^*K_E^H(1+h)}\cdot hK_E^H=\frac{bh}{2p^*(1+h)}.
  \]
  Switching to polar coordinates $(\xi,\eta) = (r\cos\theta, r\sin\theta)$, the estimated amplitude of the limit cycle at $K_E$ can be obtained from the ODE
  \begin{equation}\label{eq:rprime}
  r'=\alpha(K_E)r+l_1r^3+O(r^4)
  \end{equation}
  as
  \[
  r^*(K_E)=\sqrt{\frac{\alpha(K_E)}{|l_1|}}.
  \]
  Thus, for $h>0$,
  \begin{equation}\label{eq:rKE}
  r^*((1+h)K_E^H)=\sqrt{\frac{\frac{bh}{2p^*(1+h)}}{\frac{b^2(d-m)(p^*-1)^3}{4K_P^2p^{*}(p^*+1)^2}}}=\sqrt{\frac{2K_P^2(p^*+1)^2h}{b(1+h)(d-m)(p^*-1)^3}}=\frac{K_P(p^*+1)}{p^*-1}\sqrt{\frac{2h}{b(1+h)(d-m)(p^*-1)}}.
  \end{equation}
  This means that, near the Hopf bifurcation, the limit cycle in the original (translated) variables satisfies
  \begin{align*}
      u'(t)=& \frac{(\alpha^2+\beta^2)p^*}{a}r^*\cos\theta(t)\approx\frac{bp^*(d-m)(p^*-1)}{a(p^*+1)}r^*\cos\theta(t),\\
      v'(t)=& \alpha r^*\cos\theta(t)+\beta r^*\sin\theta(t)\approx\sqrt{\frac{b(d-m)(p^*-1)}{p^*+1}}r^*\sin\theta(t),
  \end{align*}
  and thus
  \begin{equation}\label{eq:ampli}
 A(u)\approx\frac{K_Pp^*}{a}\sqrt{\frac{2hb(d-m)}{(1+h)(p^*-1)}},\qquad A(v)\approx\frac{K_P}{p^*-1}\sqrt{\frac{2h(p^*+1)}{(1+h)}}.
  \end{equation}
  The switch to polar coordinates also allows us to obtain an estimate for the period of the oscillations. Indeed, we obtain the ODE
  \[
  \theta'=-\beta+c_2r^2+O(r^4),
  \]
  where (see Appendix \ref{sec:c2}) 
\begin{equation*}
\setlength\arraycolsep{0.1em}\begin{array}{rcl}
c_2 &=& \cfrac{1}{16}\left[(g_{\xi\eta\eta}-f_{\eta\eta\eta})
  +\cfrac{1}{3\beta}\bigl(2f_{\xi\eta}^2+2g_{\eta\eta}^2-5f_{\xi\eta}g_{\eta\eta}
  +5f_{\eta\eta}^2-f_{\eta\eta}g_{\xi\eta}+2g_{\xi\eta}^2\bigr)\right]\\[4mm]
  &=&\cfrac{1}{16}\Bigg[\left(\cfrac{\beta^3p^*}{a}G_{uvv}-\cfrac{a\beta}{p^*}F_{vvv}\right)\\[4mm]
  &&+\cfrac{1}{3\beta}\left(2\beta^2F_{uv}^2+2\beta^2G_{vv}^2-5\beta^2F_{uv}G_{vv}
   +5\left(\cfrac{a}{p^*}\right)^{2}F_{vv}^2-\beta^2F_{vv}G_{uv}
   +2\left(\cfrac{\beta^2p^*}{a}\right)^{2}G_{uv}^2\right)\Bigg]\\[4mm]
     &=&\cfrac{1}{16}\Bigg[\left(\cfrac{2\beta (p^*-1)^3}{K_P^2p^{*2}}\left(\beta^2-\cfrac{3b(d-m)}{p^*+1}\right)\right)
  +\cfrac{\beta}{3K_P^2p^{*2}}\Bigg(2(d-m)^2(p^*-1)^4
+\cfrac{8b^2(p^*-1)^2}{(p^*+1)^2}\\[4mm]
&&+\cfrac{10b(d-m)(p^*-1)^3}{p^*+1}
-\cfrac{2b(d-m)(p^*-1)^4}{p^*+1}\Bigg)+\cfrac{20\,b^2(d-m)^2(p^*-1)^4}{3\beta\,K_P^2p^{*2}(p^*+1)^2}
+\cfrac{2\beta^3(p^*-1)^4}{3K_P^2p^{*2}}\Bigg]\\[4mm]
&=&\cfrac{2\beta^3(p^*-1)^2(p^*-4)}{K_P^2p^{*2}}
+\cfrac{\beta}{3K_P^2p^{*2}}\left[2(d-m)^2(p^*-1)^4+\cfrac{8b^2(p^*-1)^2}{(p^*+1)^2}
+2\beta^2(p^*-1)^2(6-p^*)\right]\\[4mm]
&&+\cfrac{20\beta^3(p^*-1)^2}{3K_P^2p^{*2}}
+\cfrac{2\beta^3(p^*-1)^4}{3K_P^2p^{*2}}\\[4mm]
&=&\cfrac{\beta(p^*-1)^2}{K_P^2p^{*2}}\left[2\beta^2(p^*-4)+\cfrac{2}{3}(d-m)^2(p^*-1)^2+\cfrac{8b^2}{3(p^*+1)^2}
+\cfrac{2}{3}\beta^2(6-p^*)+\cfrac{20}{3}\beta^2+\cfrac{2}{3}\beta^2(p^*-1)^2\right]\\[4mm]
&=&\cfrac{\beta(p^*-1)^2}{K_P^2p^{*2}}\left[\cfrac{2}{3}(d-m)^2(p^*-1)^2+\cfrac{8b^2}{3(p^*+1)^2}
+\beta^2\left(2(p^*-4)+\cfrac{2}{3}(6-p^*)+\cfrac{20}{3}+\cfrac{2}{3}(p^*-1)^2\right)\right]\\[4mm]
&=&\cfrac{\beta(p^*-1)^2}{3K_P^2p^{*2}}\left[2(d-m)^2(p^*-1)^2+\cfrac{8b^2}{(p^*+1)^2}
+2\beta^2\Big(3(p^*-4)+(6-p^*)+10+(p^*-1)^2\Big)\right].\\[4mm]
\end{array}
\end{equation*}
From \eqref{eq:rKE}, at $K_E=(1+h)K_E^H$ we get
\[
\begin{aligned}
c_2r^{*2}(K_E)
=&\frac{\beta(p^*-1)^2}{24K_P^2p^{*2}}
\left[\beta^2(p^{*2}+5)+(d-m)^2(p^*-1)^2+\frac{4b^2}{(p^*+1)^2}\right]
\cdot\frac{2hK_P^2(p^*+1)^2}{b(d-m)(1+h)(p^*-1)^3}\\
=&\frac{\beta h(p^*+1)^2}{12p^{*2}(1+h)b(d-m)(p^*-1)}
\left[\beta^2(p^{*2}+5)+(d-m)^2(p^*-1)^2+\frac{4b^2}{(p^*+1)^2}\right]\\
=&\frac{\beta h(p^*+1)^2}{12p^{*2}(1+h)\beta^2(p^*+1)}
\left[\beta^2(p^{*2}+5)+(d-m)^2(p^*-1)^2+\frac{4b^2}{(p^*+1)^2}\right]\\
=&\frac{\beta h(p^*+1)^2}{12p^{*2}(1+h)(p^*+1)}
\left[(p^{*2}+5)+\frac{(d-m)(p^*-1)(p^*+1)}{b}+\frac{4b}{(d-m)(p^*-1)(p^*+1)}\right],
\end{aligned}
\]
and therefore the period $T$ can be estimated as
\begin{equation}\label{eq:period}
\begin{aligned}
T\approx&\frac{2\pi}{\beta-c_2r^{*2}}\\
=&\frac{2\pi}{\beta}\left(1+\frac{h\,(p^*+1)}{12\,p^{*2}(1+h)}
\left[(p^{*2}+5)+\frac{(d-m)(p^*-1)(p^*+1)}{b}+\frac{4b}{(d-m)(p^*-1)(p^*+1)}\right]\right).
\end{aligned}
\end{equation}
Observe, in particular, that the above estimate is independent of $K_P$.
\section{Numerical simulations}\label{sec:simulations}
\begin{figure}
\centering
\includegraphics[scale=0.65]{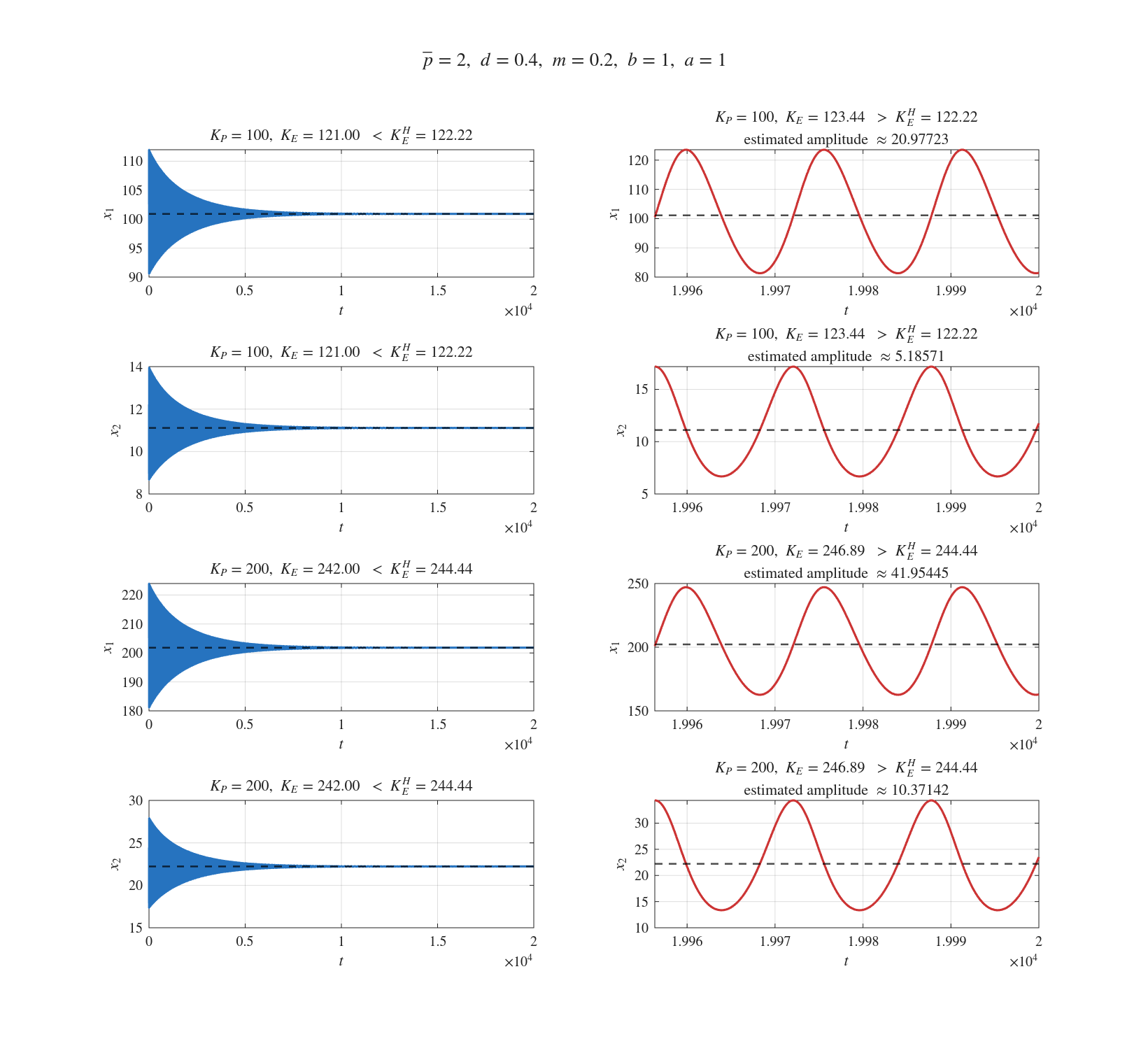}
\caption{Time-integration of \eqref{eq:Wang} from $t=0$ up to $t=10\,000$, for parameter values $\overline p=2,\,d=0.4,\,m=0.2,\,a=b=1$, as well as $K_P=100$ (left) and $K_P=200$ (right) at $K_E=0.99K_E^H$ (top 4 subfigures) and $K_E=1.01K_E^H$ (bottom 4 subfigures).}
\label{fig:hopf}
\end{figure}

\begin{figure}
\centering
\includegraphics[scale=0.75]{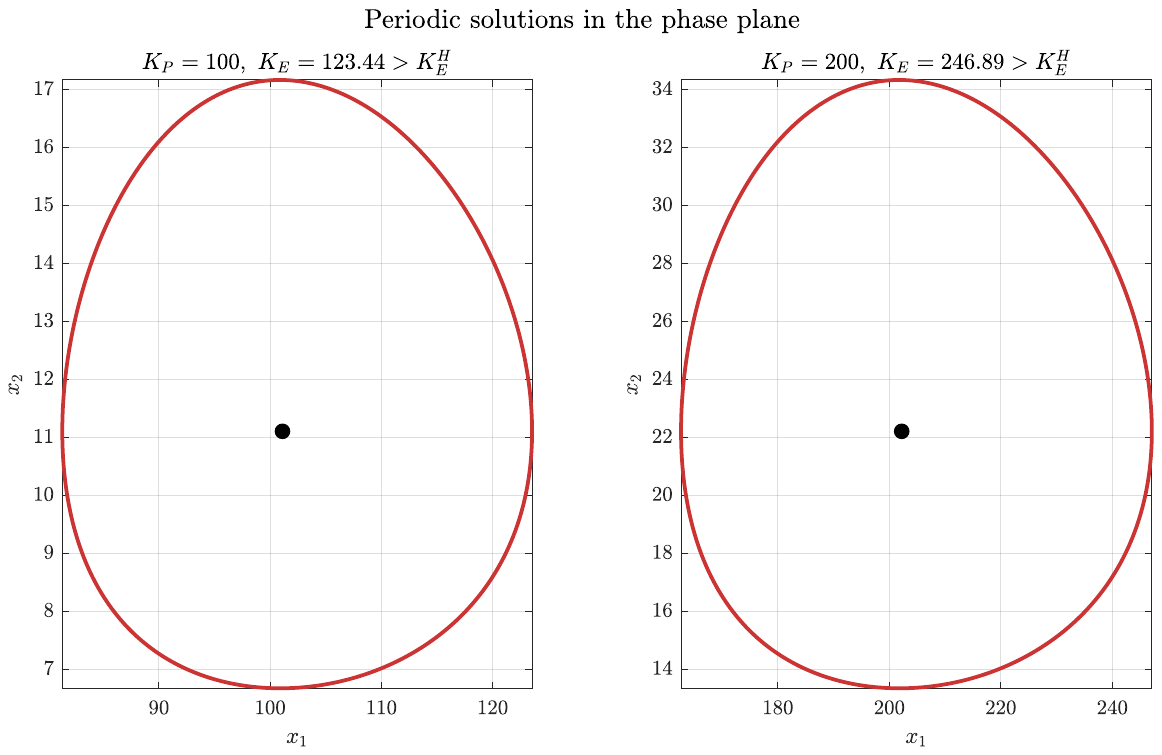}
\caption{Limit cycle of \eqref{eq:Wang} at $K_E=0.99K_E^H$ for $\overline p=2,\,d=0.4,\,m=0.2,\,a=b=1$,  as well as $K_P=100$ (left) and $K_P=200$ (right).}
\label{fig:periodicx1x2}
\end{figure}

 The value $K_E^H$, as well as the estimated amplitudes and periods of the limit cycles for slightly larger values of $K_E$, can also be confirmed by means of numerical simulations. In particular, Figure \ref{fig:hopf} shows the result of integrating \eqref{eq:Wang} through the MATLAB routine \texttt{ode45}. As in \citep{Wang2015}, we consider the parameter values $\overline p=2,\,d=0.4,\,m=0.2,\,a=b=1$. As for the parameter $K_P$, we consider the value $100$ (four top subfigures) and $200$ (four bottom subfigures). For both values, we show that the system with $K_E=0.99K_E^H$ converges to the equilibrium (left subfigures) and that with $K_E=1.01K_E^H$ converges to a limit cycle (right subfigures). In the latter case, we remark that the amplitude obtained in the simulation is very close to the estimate \eqref{eq:ampli} obtained by means of the Hopf normal form. The limit cycle in the phase plane, for both values of $K_P$, is shown in Figure \ref{fig:periodicx1x2}.

\begin{figure}
\centering
\includegraphics[scale=0.9]{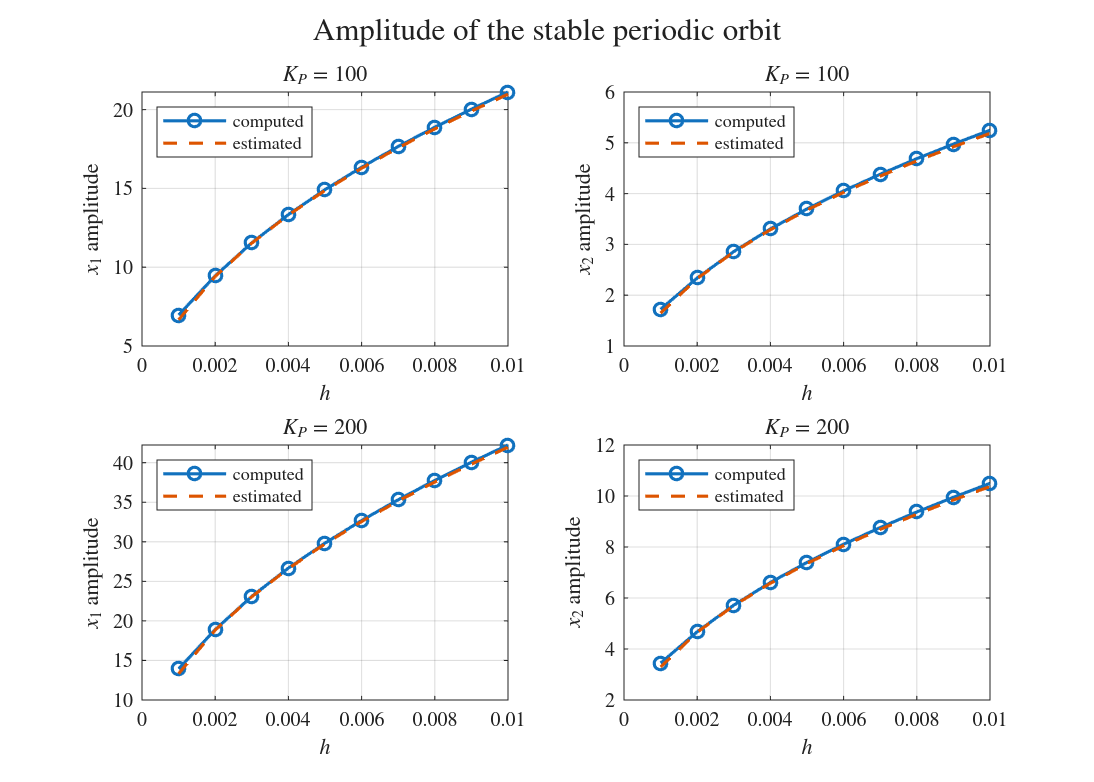}
\caption{Comparison between estimated amplitudes \eqref{eq:ampli} of periodic solutions of \eqref{eq:Wang} at $(1+h)K_E^H$ for $h=k\cdot 10^{-3}$, $k=1,\ldots,10$, for $K_P=100$ (left) and $K_P =200$ (right).}
\label{fig:amplitude}
\end{figure}

\begin{figure}
\centering
\includegraphics[scale=0.9]{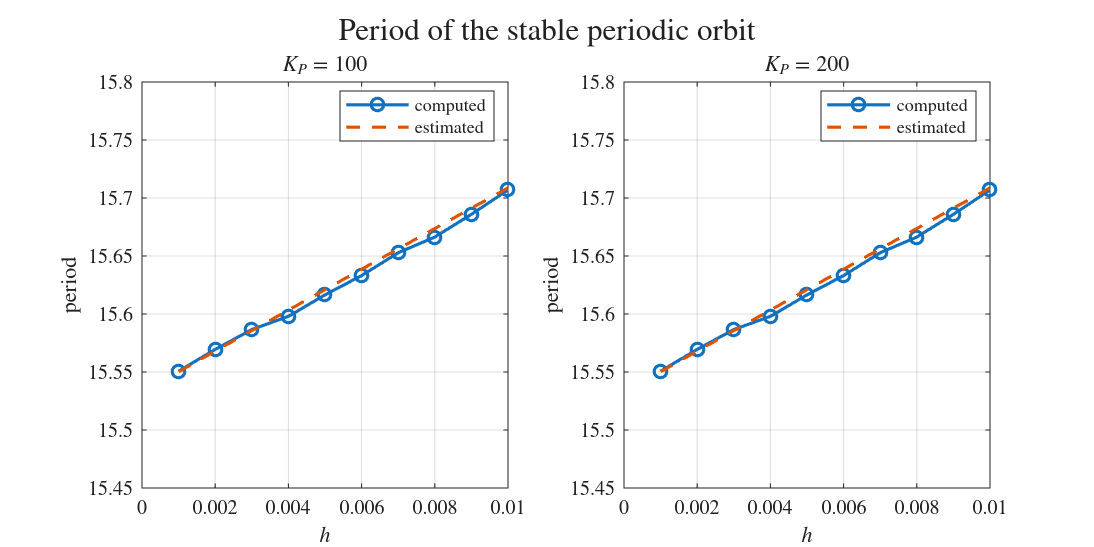}
\caption{Comparison between estimated periods \eqref{eq:period} of periodic solutions of \eqref{eq:Wang} at $(1+h)K_E^H$ for $h=k\cdot 10^{-3}$, $k=1,\ldots,10$, for $K_P=100$ (left) and $K_P =200$ (right).}
\label{fig:period}
\end{figure}

The next numerical experiments show a more systematic comparison between the amplitudes and periods estimated via the Hopf normal form and those computed by means of numerical time-integration. We performed both experiments integrating up to $t=20\,000$, at which the profile of the solution visually appears approximately periodic for all tested values of $K_E$ and $K_P$. For the former, we choose the values $(1+h)K_E^H$ for $h=k\cdot 10^{-3}$, $k=1,\ldots,10$. As for the latter, we consider again the cases $K_P=100$ and $K_P =200.$ Figure \ref{fig:amplitude} shows the comparison between the estimated amplitudes of both components, defined as in \eqref{eq:ampli}, and the corresponding computed amplitudes. Since the amplitude of each component of the solution increases over time while converging to the limit cycle, the computed amplitudes are defined as the differences between the maximum and the minimum values of the solution through the entire time series. Finally, Figure \ref{fig:period}, shows the comparison between the estimated period, defined as in \eqref{eq:period}, and the computed period. The latter has been obtained as the average of the last 10 differences in the time series between a local minimum and the previous. Local minima have been computed by means of the MATLAB routine \texttt{islocalmin}. Observe that the estimated period is independent of the chosen value of $K_P$, as also remarked in Section \ref{sec:hopf}.
\section{Conclusions}\label{sec:conclusion}
In this paper we revisited the two-dimensional product–resource supply chain system obtained from the reduction of the model by \cite{Wang2015} and showed that its equilibrium can lose stability through a Hopf bifurcation as the environment carrying capacity crosses a critical value, which depends on other model parameters. We computed this value, and established the direction and stability of the resulting periodic orbit via the first Lyapunov coefficient, closing a gap in the original analysis. In particular, we found that increasing the environmental capacity or the maximum production rate destabilizes the equilibrium and pushes the system toward oscillation, while increasing the resource level at which production reaches half its maximum rate, or increasing the remanufacturing rate, restores stability.

Whether the periodic regime should be read as a failure of the system to reach its intended operating point, or as a benign and expected feature, depends on what the model's product variable represents. Normally speaking, sustained oscillation in production and resource use is undesirable for the reasons mentioned in the Introduction. However, for a product with inherently seasonal demand (such as beachwear or holiday decorations), it may instead indicate a rhythm that the production is actually expected to follow \citep{Soysal2012}.

The Hopf bifurcation established in the present paper originates exclusively from the coupling between the product and resource variables, with no delay involved. In practice, however, due to manufacturing, procurement, or processing times, production in the present depends rather on the resource conditions at some past time, and this mechanism is the most commonly used in the supply chain literature to explain oscillatory dynamics, as explained in the Introduction. In the future, we plan to extend the current work by introducing a delay in the effective production rate term, and analyze how the two mechanisms for generating periodic solutions interact with each other. In particular, it would be interesting to investigate whether the delays narrows (as we expect) or widens the parameter region in which the nontrivial equilibrium is stable, and whether it introduces entirely new dynamics, such as a second Hopf bifurcation. Finally, the coupling of these mechanisms with possible stochastic effects in the model, introduced by \cite{DAmbrosio2026}, will also be worth exploring.
\section*{Acknowledgements}
The authors are members of INdAM research group GNCS and of UMI research group “Modellistica socio-epidemiologica”. The work was partially supported by the Italian Ministry of University and Research (MUR) through the PRIN 2022 project (No. 20229P2HEA) “Stochastic numerical modelling for sustainable innovation”, Unit of Udine (CUP G53C24000710006).
\appendix
\section{Derivation of $c_2$}\label{sec:c2}
Recalling the switch to polar coordinates $\xi = r\cos\theta$, $\eta = r\sin\theta$, at $\alpha=0$ we obtain from \eqref{eq:uv3} the equation
\begin{equation}\label{eq:rthetatrue}\begin{aligned}
   r' &= \chi_2(\theta)\,r^2+\chi_3(\theta)\,r^3+O(r^4),\\
  \theta' &= -\beta+\psi_2(\theta)\,r+\psi_3(\theta)\,r^2+O(r^3),
  \end{aligned}
  \end{equation}
where 
\begin{align*}
\chi_2(\theta) &= f_{\xi\eta}\cos^2\theta\sin\theta+\left(\frac12f_{\eta\eta}+g_{\xi\eta}\right)\sin^2\theta\cos\theta
 +\frac12g_{\eta\eta}\sin^3\theta,\\
\psi_2(\theta) &= g_{\xi\eta}\cos^2\theta\sin\theta+\left(\frac12g_{\eta\eta}-f_{\xi\eta}\right)\sin^2\theta\cos\theta
  -\frac12f_{\eta\eta}\sin^3\theta,
\end{align*}
and
\begin{align*}
\chi_3(\theta) &= \frac12f_{\xi\eta\eta}\cos^2\theta\sin^2\theta+\left(\frac16f_{\eta\eta\eta}+\frac12g_{\xi\eta\eta}\right)\sin^3\theta\cos\theta
  +\frac16g_{\eta\eta\eta}\sin^4\theta,\\
\psi_3(\theta) &= \frac12g_{\xi\eta\eta}\cos^2\theta\sin^2\theta+\left(\frac16g_{\eta\eta\eta} -\frac12f_{\xi\eta\eta}\right)\sin^3\theta\cos\theta
 -\frac16f_{\eta\eta\eta}\sin^4\theta.
\end{align*}

Let $\langle\cdot\rangle$ denote the average of a function in the interval $[0,2\pi]$, defined as $\frac1{2\pi}\int_0^{2\pi}(\cdot)\,\mathrm{d}\theta$. Observe that $\langle\chi_2\rangle=\langle\psi_2\rangle=0$, $\langle\chi_3\rangle=\frac1{16}(f_{\xi\eta\eta}+g_{\eta\eta\eta})$ and $\langle\psi_3\rangle=\frac{1}{16}(g_{\xi\eta\eta}-f_{\eta\eta\eta})$. The former justifies the absence of a quadratic term in \eqref{eq:rprime}: the idea is that near the bifurcation $r$ is small and changes slowly (on timescale $1/|\alpha|$ or $1/|l_1|r^2$), while $\theta$ oscillates rapidly at frequency $\beta$ (since $\theta'\approx -\beta$). \eqref{eq:rprime} is obtained by replacing the true right-hand side of $r'$ in \eqref{eq:rthetatrue} by its mean over one fast period $2\pi/\beta$.

We now introduce new coordinates $(\rho,\vartheta)$ which satisfy
\begin{equation}\label{eq:transform}
  r = \rho + \mu(\vartheta)\,\rho^2,\qquad
  \theta = \vartheta + \nu(\vartheta)\,\rho
\end{equation}
for suitable $2\pi$-periodic function $\mu(\vartheta)$ and $\nu(\vartheta)$, such that $\rho'$ has no quadratic term in $r$ and $\theta'$ has no linear term in $r$. Let us write $\rho' =B(\vartheta)\rho^3+O(\rho^4)$ and
$\vartheta' = -\beta+C(\vartheta)\rho^2+O(\rho^3)$, for unknowns $B,\,C$ to be found. Differentiating \eqref{eq:transform} and ignoring the higher-order terms gives 
\[
  r' = \rho'+\mu'(\vartheta)\vartheta'\,\rho^2+2\mu(\vartheta)\rho\rho',\qquad
  \theta' = \vartheta'+\nu'(\vartheta)\vartheta'\,\rho+\nu(\vartheta)\rho'.
\]
Using the expressions for $\rho'$ and $\vartheta$, we get
\begin{equation}\label{eq:difftransform}
\begin{aligned}
r' &= -\beta \mu'(\vartheta)\rho^2+B(\vartheta)\rho^3+O(\rho^4),\\
\theta' &= -\beta-\beta \nu'(\vartheta)\rho+C(\vartheta)\rho^2+O(\rho^3).
\end{aligned}
\end{equation}
Writing the right-hand sides of \eqref{eq:rthetatrue} in terms of the new coordinates gives
\begin{equation}
\begin{aligned}\label{eq:rthetatruetransform}
r' &= \chi_2(\vartheta)\rho^2+[2\chi_2(\vartheta)\mu(\vartheta)+\chi_2'(\vartheta)\nu(\vartheta)+\chi_3(\vartheta)]\rho^3(\vartheta)+O(\rho^4),\\
\theta' &= -\beta+\psi_2(\vartheta)\rho+[\psi_2(\vartheta)\mu(\vartheta)+\psi_2'(\vartheta)\nu+\psi_3(\vartheta)]\rho^2+O(\rho^3).
\end{aligned}
\end{equation}
Comparing \eqref{eq:difftransform} and \eqref{eq:rthetatruetransform}, we obtain $\mu'=-\frac{\chi_2}{\beta}$ and $\nu'=-\frac{\psi_2}{\beta}$. Consequently,
\begin{equation}\label{eq:D3E2}
  B = 2\chi_2 \mu+\chi_2'\nu+\chi_3,\qquad
  C = \psi_2 \mu + \psi_2'\nu+\psi_3.
\end{equation}
Averaging \eqref{eq:D3E2} over $\vartheta$ gives the normal form
$\dot\rho=l_1\rho^3+O(\rho^4)$, $\dot\vartheta=-\beta+c_2\rho^2+O(\rho^3)$
with
\begin{equation*}
l_1=2\langle\chi_2\mu\rangle+\langle\chi_2'\nu\rangle+\langle\chi_3\rangle,
  \qquad
c_2=\langle\psi_2\mu\rangle+\langle\psi_2'\nu\rangle+\langle\psi_3\rangle.
\end{equation*}
The expression on the left corresponds to the first line of \eqref{eq:lyap1}, as obtained explicitly in \citep{Guckenheimer2013}. We will now follow similar arguments in order to expand the expression on the right. Integration by parts gives
\[
  \langle \psi_2' \nu\rangle = -\langle\psi_2\nu'\rangle
  = 0-\Big\langle\psi_2\cdot\Big(-\frac{\psi_2}\beta\Big)\Big\rangle
  = \frac1\beta\langle\psi_2^2\rangle,
\]
where the first equality follows from $2\pi$-periodicity of $\psi_2$ and $\nu$. In order to compute $\langle\psi_2\mu\rangle$ and $\langle\psi_2^2\rangle$, let us first observe that, for $m,\,n\in\mathbb{N}$, we have 
\begin{equation}\label{eq:zero-avg}
  \langle\cos^m\theta\sin^n\theta\rangle = 0 \quad\text{if $m$ or $n$ is odd.}
\end{equation}
Indeed, if $n=2k+1$, then
\[
\int_0^{2\pi}\cos^m\theta\sin^{2k+1}\theta\mathrm{d}\theta=\int_0^{2\pi}\cos^m\theta(1-\cos^2\theta)^k\sin\theta\mathrm{d}\theta=-\int_1^1x^m(1-x^2)^k\mathrm{d}x=0,
\]
and similarly for the case where $m$ is odd. Moreover, $\langle\sin^2\theta\rangle=\langle\cos^2\theta\rangle=\tfrac12$,
$\langle\sin^4\theta\rangle=\langle\cos^4\theta\rangle=\tfrac38$,
$\langle\sin^6\theta\rangle=\langle\cos^6\theta\rangle=\tfrac5{16}$,
$\langle\sin^2\theta\cos^2\theta\rangle=\tfrac18$, and $\langle\sin^4\theta\cos^2\theta\rangle=\langle\cos^4\theta\sin^2\theta\rangle=\frac1{16}$. In order to obtain the expression of $\mu$, observe that
\[
\int\cos^2\theta\sin\theta\,\mathrm{d}\theta=-\tfrac13\cos^3\theta,\quad\int\cos\theta\sin^2\theta\,\mathrm{d}\theta=\tfrac13\sin^3\theta,\quad\int\sin^3\theta\,\mathrm{d}\theta
  =-\cos\theta+\tfrac13\cos^3\theta.
\]
Recalling that $\mu=-\frac{\chi_2}{\beta}$, this implies
\[
\mu(\theta) = \frac1\beta\left[\left(\frac13f_{\xi\eta}-\frac16g_{\eta\eta}\right)\cos^3\theta+\frac12g_{\eta\eta}\cos\theta-\left(\frac16f_{\eta\eta}+\frac13g_{\xi\eta}\right)\sin^3\theta\right].
\]
By \eqref{eq:zero-avg}, this means
\[
\begin{aligned}
\langle\psi_2\mu\rangle=&\frac1{\beta}\Bigg[\left(\frac13f_{\xi\eta}-\frac16g_{\eta\eta}\right)\left(\frac12g_{\eta\eta}-f_{\xi\eta}\right)\langle\sin^2\cos^4\rangle+\frac12g_{\eta\eta}\left(\frac12g_{\eta\eta}-f_{\xi\eta}\right)\langle\sin^2\cos^2\rangle\\
-&g_{\xi\eta}\left(\frac16f_{\eta\eta}+\frac13g_{\xi\eta}\right)\langle\sin^4\cos^2\rangle+\frac12f_{\eta\eta}\left(\frac16f_{\eta\eta}+\frac13g_{\xi\eta}\right)\langle\sin^6\rangle\Bigg]\\
=&\frac1{\beta}\Bigg[\frac1{16}\left(-\frac13f_{\xi\eta}^2-\frac1{12}g_{\eta\eta}^2+\frac13f_{\xi\eta}g_{\eta\eta}\right)+\frac1{16}\left(\frac1{2}g_{\eta\eta^2}-f_{\xi\eta}g_{\eta\eta}\right)-\frac1{16}\left(\frac16f_{\eta\eta}g_{\xi\eta}+\frac13g_{\xi\eta}^2\right)\\
+&\frac5{32}\left(\frac16f_{\eta\eta}^2+\frac13f_{\eta\eta}g_{\xi\eta}\right)\Bigg]\\
=&\frac1{192\beta}\left[-4f_{\xi\eta}^2+5g_{\eta\eta}^2-8f_{\xi\eta}g_{\eta\eta}+5f_{\eta\eta}^2+8f_{\eta\eta}g_{\xi\eta}-4g_{\xi\eta}^2\right]
\end{aligned}
\]
and
\[
\begin{aligned}
 \langle\psi_2^2\rangle=&\left[g_{\xi\eta}^2\langle\sin^2\cos^4\rangle+\left(\frac12g_{\eta\eta}-f_{\xi\eta}\right)^2\langle\sin^4\cos^2\rangle+\frac14f_{\eta\eta}^2\langle\sin^6\rangle-f_{\eta\eta}g_{\xi\eta}\langle\sin^4\cos^2\rangle\right]\\
 =&\left[\frac1{16}\left(g_{\xi\eta}^2+\frac14g_{\eta\eta}^2+f_{\xi\eta}^2-f_{\xi\eta}g_{\eta\eta}-f_{\eta\eta}g_{\xi\eta}\right)+\frac5{64}f_{\eta\eta}^2\right]\\
 =&\frac1{64}\left[4g_{\xi\eta}^2+g_{\eta\eta}^2+4f_{\xi\eta}^2-4f_{\xi\eta}g_{\eta\eta}-4f_{\eta\eta}g_{\xi\eta}+5f_{\eta\eta}^2\right].
\end{aligned}
\]
Eventually, we get
\[
\begin{aligned}
c_2=&\frac1{16}(f_{\xi\eta\eta}+g_{\eta\eta\eta})+\langle\psi_2\mu\rangle+\frac1{\beta}\langle\psi_2^2\rangle\\
=&\frac1{16}(f_{\xi\eta\eta}+g_{\eta\eta\eta})\\
+&\frac1{192\beta}\left[-4f_{\xi\eta}^2+5g_{\eta\eta}^2-8f_{\xi\eta}g_{\eta\eta}+5f_{\eta\eta}^2+8f_{\eta\eta}g_{\xi\eta}-4g_{\xi\eta}^2+3(4g_{\xi\eta}^2+g_{\eta\eta}^2+4f_{\xi\eta}^2-4f_{\xi\eta}g_{\eta\eta}-4f_{\eta\eta}g_{\xi\eta}+5f_{\eta\eta}^2)\right]\\
=&\frac1{16}(f_{\xi\eta\eta}+g_{\eta\eta\eta})+\frac1{192\beta}\left[8f_{\xi\eta}^2+8g_{\eta\eta}^2-20f_{\xi\eta}g_{\eta\eta}+20f_{\eta\eta}^2-4f_{\eta\eta}g_{\xi\eta}+8g_{\xi\eta}^2\right]\\
=&\frac1{16}\left[(g_{\xi\eta\eta}-f_{\eta\eta\eta})+\frac1{3\beta}(2f_{\xi\eta}^2+2g_{\eta\eta}^2-5f_{\xi\eta}g_{\eta\eta}+5f_{\eta\eta}^2-f_{\eta\eta}g_{\xi\eta}+2g_{\xi\eta}^2)\right].
\end{aligned}
\]
\bibliographystyle{plainnat}
\bibliography{supply}

\end{document}